\documentclass{article}
\usepackage{amsmath,amsfonts,amssymb,amsbsy}
\usepackage[pdftex]{graphicx}
\begin{document} 

\centerline{\large\bf 
Asymptotic modeling of the JKR adhesion contact}
\medskip
\centerline{\large\bf  for a thin elastic layer}
\medskip
\medskip
\centerline{I.\,I.~ARGATOV, G.\,S.~MISHURIS\footnote{Corresponding author, Email: ggm@aber.ac.uk} and V.\,L.~POPOV}
\medskip
\medskip
\centerline{\it
Institute of Mathematics and Physics, Aberystwyth University, SY23 3BZ, Wales, UK}
\centerline{\it
Institut f{\"u}r Mechanik, Technische Universit{\"a}t Berlin, 10623 Berlin, Germany}

\bigskip
\medskip
{\bf Abstract:}
The Johnson--Kendall--Roberts model of frictionless adhesive contact is extended to the case of a thin transversely isotropic elastic layer bonded to a rigid base. The leading-order asymptotic models are obtained for both compressible and incompressible elastic layers. The boundary conditions for the contact pressure approximation on the contour of the contact area have been derived from the boundary-layer solutions, which satisfy the condition that the stress intensity factor of the contact pressure density should have the same value all round the contact area boundary. In the incompressible case, a perturbation solution is obtained for a slightly perturbed circular contact area.

\medskip

\setcounter{equation}{0}

\section{Introduction}
\label{SecI}

Mechanical aspects of adhesion in biological systems have been received much attention in previous years \cite{SpolenakGorbGao2005,FilippovPopovGorb2011}.
In particular, it was found that adhesion of live cells to external surfaces plays an important role in many cellular processes \cite{SafranGovNicolas2005,KendallKendallRehnfeldt2011}. Recently, the effect of contact surface shape as a factor to control focal adhesion lifetime was studied \cite{RizzaQianGao2011}, while modeling the cell and substrate as a semi-infinite elastic media.

In the present paper, we consider the classical JKR (Johnson--Kendall--Roberts) model \cite{JohnsonKendallRoberts1971}, which was originally developed for a spherical rigid punch in frictionless adhesive contact with an isotropic elastic half-space, and extend it to the case of a thin transversely isotropic elastic layer. Previously, asymptotic models for axisymmetric adhesive indentation of a thin isotropic elastic layer were developed in the 
incompressible \cite{Yang2002} and compressible \cite{Reedy2006,Yang2006} cases.
A numerical analysis of the adhesive contact between a soft elastic layer and a rough rigid
substrate was developed in \cite{CarboneMangialardi2008} by using a Green's function approach.

The main results of the present paper are represented by the leading-order asymptotic models of non-axisymmetric adhesive contact for a thin elastic layer bonded to a rigid base. In the case of a thin {\it compressible\/} layer indented by a rigid punch with the shape function $\varphi({\bf y})$, the contact pressure, $p({\bf y})$, which is distributed over the contact area $\omega$ bounded by the contour $\Gamma$, is determined by the equations
$$
p({\bf y})=\frac{A_{33}}{h}\bigl(\delta_0-\varphi({\bf y})\bigr),\quad {\bf y}\in\omega,
$$
\begin{equation}
p({\bf y})=-\sqrt{\frac{2A_{33}\Delta\gamma}{h}},\quad {\bf y}\in\Gamma.
\label{jkr1(1.1)}
\end{equation}
Here, $\delta_0$ is the punch displacement, $A_{33}$ is the aggregate elastic modulus of the elastic layer, $h$ is the layer thickness, $\Delta\gamma$ is the surface energy. 

In the case of a thin bonded {\it incompressible\/} layer, the following boundary-value problem is obtained for the contact pressure approximation:
$$
-\frac{h^3}{3G^\prime}\Delta_y p({\bf y})
=\delta_0-\varphi({\bf y}),\quad {\bf y}\in\omega,
$$
\begin{equation}
p({\bf y})=0,\quad 
\frac{\partial p}{\partial n}({\bf y})
=-\frac{1}{h^{3/2}}\sqrt{6G^\prime\Delta\gamma},\quad {\bf y}\in\Gamma,
\label{jkr1(1.2)}
\end{equation}
Here, $G^\prime$ is the out-of-plane shear elastic modulus of the elastic layer, $\Delta_y$ is the Laplace differential operator, $\partial/\partial n$ is the inward normal derivative.

It is interesting that the boundary conditions (\ref{jkr1(1.1)}) and (\ref{jkr1(1.2)}) have been derived from the boundary-layer solutions, which satisfy the condition \cite{JohnsonGreenwood2005} that the stress intensity factor of the contact pressure density should have the same value all round the contact area boundary.

The rest of the paper is organized as follows. In Section~\ref{Sec2}, we formulate the  unilateral contact problem of frictionless adhesive contact with a transversely isotropic bonded to a rigid base. The case of a thin compressible layer is considered in Section~\ref{Sec3}, and the elliptic frictionless adhesive contact is solved in detail. In Section~\ref{Sec4}, we study the case of a thin incompressible layer, and a perturbation solution is obtained for a slightly perturbed circular contact area. In dealing with the perturbed unilateral adhesive contact problem, we utilize the developed previously asymptotic technique \cite{ArgatovMishuris2011p}. Finally, in Section~\ref{SecDC}, we formulate our conclusions. 

\section{Adhesion contact problem formulation for an elastic layer}
\label{Sec2}

We consider the frictionless adhesive contact between a transversely isotropic elastic layer of a relatively small thickness, $h$, bonded to a rigid base and a rigid punch described by the shape function
$$
z=\varphi(y_1,y_2).
$$
As a special case we consider the punch in the form of an elliptic paraboloid 
\begin{equation}
\varphi({\bf y})=\frac{y_1^2}{2R_1}+\frac{y_2^2}{2R_2},
\label{jkr1(2.1)}
\end{equation}
where $R_1$ and $R_2$ are the radii of curvature of the principal normal cross-sections of the punch's surface at its vertex.

Let $\delta_0$ and $p({\bf y})$ denote the punch indentation depth and the contact pressure density, respectively. By making use of the standard two-dimensional Fourier transform technique, the contact problem can be reduced to the following governing integral equation \cite{Vorovich_et_al_1974,Sneddon1995}:
\begin{equation}
\frac{1}{\pi^2 h\theta}\iint\limits_{\omega}
p({\bf y}^\prime)K(y_1-y_1^\prime,y_2-y_2^\prime)\,d{\bf y}^\prime
=\delta_0-\varphi({\bf y}).
\label{jkr1(2.2)}
\end{equation}
Here the kernel is given by the integral 
\begin{equation}
K({\bf y})=\iint\limits_0^{{}\ \ \ +\infty}
\frac{\mathcal{L}(\lambda)}{\lambda}\cos\frac{\lambda_1 y_1}{h}\cos\frac{\lambda_2 y_2}{h}\,d\lambda_1 d\lambda_2.
\label{jkr1(2.2a)}
\end{equation}
Note that we follow the same notation as in the book \cite{ArgatovMishuris2015}, wherever possible.

In the case of a transversely isotropic elastic layer bonded to a flat rigid base, the elastic constant $\theta$ and the kernel function $\mathcal{L}(\lambda)$ are given by the known solution \cite{Fabrikant2006} (see also \cite{ArgatovSabina2013} and \cite{ArgatovMishuris2015}, Section~2.1.3).

Let $\Gamma$ denote the unknown boundary of the contact area $\omega$. Introducing a natural parametrization of the contour $\Gamma$ by formulas $y_1=f_1(s)$, $y_2=f_2(s)$, where $s$ is the arc length, we will assume that when traveling along $\Gamma$ in the direction of increasing $s$-coordinate, the region $\omega$ enclosed by $\Gamma$ remains on the left. Correspondingly, the unit vector of the inward (with respect to the domain $\omega$) normal to the contour $\Gamma$ is given by
$$
{\bf n}(s)=-f_2^\prime(s){\bf e}_1+f_1^\prime(s){\bf e}_2,
$$
where the prime denotes differentiation with respect to $s$.

Then, the stress intensity factor (SIF), $K_I(s)$, at the point ${\bf y}(s)=\bigl(f_1(s),f_2(s)\bigr)$ on the contour $\Gamma$ can be evaluated as 
\begin{equation}
K_I(s)=-\lim_{n\to 0}p\bigl({\bf y}(s)+n{\bf n}(s)\bigr)\sqrt{2\pi n},
\label{jkr1(2.6)}
\end{equation}
so that $\bigl(f_1(s)+n n_1(s), y_2(s)=f_2(s)+n n_2(s)\bigr)$ are the Cartesian coordinates of the point of observation ${\bf y}(s)+n{\bf n}(s)$, which tends to the contour $\Gamma$ along the normal at the point with the natural parameter $s$.

In the case of the JKR-type adhesive contact according to Griffith's energy balance (see, e.g., \cite{Maugis1995}), the boundary $\Gamma$ of the contact area $\omega$ is determined by the following condition \cite{JohnsonGreenwood2005,Barber2013}:
\begin{equation}
K_I(s)=2\sqrt{\theta\Delta\gamma}.
\label{jkr1(2.7)}
\end{equation}
Here, $\Delta\gamma$ denotes the work of adhesion (which is taken to be the surface energy).
We note that $\theta=E^*/2$, where $E^*$ is the so-called effective elastic modulus of the layer, which in the isotropic case is equal to $E/(1-\nu^2)$ with $E$ and $\nu$ being Young's modulus and Poisson's ratio, respectively.

Finally, the contact force, $F$, is determined by the equilibrium equation 
\begin{equation}
\iint\limits_{\omega}p({\bf y})\,d{\bf y}=F.
\label{jkr1(2.8)}
\end{equation}

Thus, Eqs.~(\ref{jkr1(2.2)}) and (\ref{jkr1(2.7)}) constitute the JKR adhesive contact problem for an elastic layer, which is based on the assumption \cite{JohnsonGreenwood2005} that the SIF of the contact pressure should have the same value all round the boundary $\Gamma$.

\section{Adhesive contact in the case of a thin compressible layer}
\label{Sec3}

\subsection{Leading-order asymptotic model for a thin compressible elastic layer bonded to a rigid base}
\label{Sec3.1}

In what follows, we will make use of the expansion 
\begin{equation}
\frac{\mathcal{L}(\lambda)}{\theta \lambda}= 
\mathcal{M}_0+\mathcal{M}_1 \lambda^2+\mathcal{M}_2 \lambda^4+\ldots
\label{jkr1(2.3)}
\end{equation}
with the first coefficients given by 
\begin{equation}
\mathcal{M}_0=\frac{\mathcal{A}}{\theta}=\frac{1}{A_{33}}, \quad
\mathcal{M}_1=\frac{A_{13}(A_{13}-A_{44})}{3A_{33}^2A_{44}}.
\label{jkr1(2.4)}
\end{equation}
Note also that the dimensional asymptotic contact $\mathcal{A}=\theta\mathcal{M}_0$ is given by
\begin{equation}
\mathcal{A}=\frac{A_{11}A_{33}-A_{13}^2}{(\gamma_1+\gamma_2)A_{11}A_{33}},
\label{jkr1(2.4A)}
\end{equation}
where  $\gamma_1$ and $\gamma_2$ are the roots of the characteristic equation.

Using the distributional asymptotic analysis \cite{EstradaKanwal1990}, the following formal asymptotic expansion for the integral operator on the left-hand side of Eq.~(\ref{jkr1(2.2)}) was established \cite{Argatov2005Acta}:
\begin{equation}
\frac{1}{\pi^2 h\theta}\iint\limits_{\omega}
p({\bf y}^\prime)K({\bf y}-{\bf y}^\prime)\,d{\bf y}^\prime
\simeq 
h\sum_{n=0}^\infty(-1)^n \mathcal{M}_n h^{2n} \Delta_y^n p({\bf y}).
\label{jkr1(2.5)}
\end{equation}
Here, $\Delta_y=\partial^2/\partial y_1^2+\partial^2/\partial y_2^2$ is the Laplace operator.

It is to emphasize that the asymptotic expansion (\ref{jkr1(2.5)}) is valid in the inner region of the contact area $\omega$, that is relatively far from the boundary $\Gamma$.

Keeping in the series (\ref{jkr1(2.5)}) only the first asymptotic term and substituting the result into Eq.~(\ref{jkr1(2.2)}), we obtain 
\begin{equation}
p({\bf y})=k\bigl(\delta_0-\varphi({\bf y})\bigr),\quad {\bf y}\in\omega,
\label{jkr1(3.1)}
\end{equation}
where we have introduced the notation
$$
k=\frac{A_{33}}{h}.
$$

It is clear that Eq.~(\ref{jkr1(3.1)}) should be supplemented with some boundary condition to determine the contact area $\omega$.

\subsection{Leading-order boundary-layer solution in the compressible case}
\label{Sec3.2}

Assuming that the layer thickness $h$ is relatively small with respect to the characteristic length of the contact area $\omega$, we introduce a small dimensionless parameter, $\varepsilon$, and set 
\begin{equation}
h=\varepsilon h_*,
\label{jkr1(3.2)}
\end{equation}
where $h_*$ is independent of $\varepsilon$.

Moreover, following \cite{Argatov2005Acta}, we will make use of the normalization 
\begin{equation}
\delta_0=\varepsilon\delta_0^*,\quad R_1=\varepsilon^{-1} R_1^*,\quad R_2=\varepsilon^{-1} R_2^*,
\label{jkr1(3.2b)}
\end{equation}
where $\delta_0^*$, $R_1^*$, and $R_2^*$ are comparable with $h_*$, all being independent of $\varepsilon$. 

Further, we introduce dimensionless variables
\begin{equation}
\mbox{\boldmath $\eta$}=(\eta_1,\eta_2),\quad \eta_i=h_*^{-1}y_i, \quad i=1,2.
\label{jkr1(3.3)}
\end{equation}

Therefore, in light of (\ref{jkr1(3.2)}) and (\ref{jkr1(3.3)}), Eq.~(\ref{jkr1(2.2)}) takes the form
\begin{equation}
\iint\limits_{\omega_*}
p_*(\mbox{\boldmath $\eta$}^\prime)k\bigl(\varepsilon^{-1}
(\mbox{\boldmath $\eta$}-\mbox{\boldmath $\eta$}^\prime)\bigr)\,d\mbox{\boldmath $\eta$}^\prime
=\varepsilon^2\frac{\pi^2 \theta}{h_*}
\bigl(\delta_0^*-\varphi^*(\mbox{\boldmath $\eta$})\bigr),
\label{jkr1(3.4)}
\end{equation}
where $\omega_*$ is the contact area in the dimensionless coordinates (\ref{jkr1(3.3)}), and in the special case (\ref{jkr1(2.1)}) we have
$\varphi^*(\mbox{\boldmath $\eta$})=h_*^2\bigl((2R_1^*)^{-1}\eta_1^2
+(2R_2^*)^{-1}\eta_2^2\bigr)$, while  
$\varphi^*(\mbox{\boldmath $\eta$})=\varepsilon^{-1}\varphi(h_*\mbox{\boldmath $\eta$})$ in the general case.

The kernel $k(\mbox{\boldmath $\xi$})$ in Eq.~(\ref{jkr1(3.4)}) is obtained from the kernel $K({\bf y})$ (see Eqs.~(\ref{jkr1(2.2)}), (\ref{jkr1(2.2a)})) in the form
$$
k(\mbox{\boldmath $\xi$})=\frac{1}{4}\iint\limits_{-\infty}^{{}\ \ \ +\infty}
\frac{\mathcal{L}(s)}{s}e^{{\rm i}(s_1\xi_1+s_2\xi_2)} ds_1 ds_2.
$$

Finally, we introduce the so-called ``fast'' normal coordinate
\begin{equation}
\nu=\frac{n}{h}=\varepsilon^{-1}\frac{n}{h_*},
\label{jkr1(3.5)}
\end{equation}
keeping the scale of the dimensionless coordinate $s_*=s/h_*$ along $\Gamma_*$ unchanged. 

Following the asymptotic procedure developed in 
\cite{ArgatovMishuris2015,Argatov2005Acta,Alexandrov1969-2}, we arrive at the following boundary-layer integral equation:
\begin{equation}
\int\limits_0^{+\infty}
q^*_\varepsilon(s_*,\nu^\prime)M(\nu^\prime-\nu)\,d\nu^\prime=\frac{\pi\theta}{h_*}
\bigl[\delta_0^*-\varphi^*\bigl(\mbox{\boldmath $\eta$}(s_*)\bigr)\bigr].
\label{jkr1(3.6)}
\end{equation}
Here, $q^*_\varepsilon(s_*,\nu)$ is the leading-order asymptotic approximation for the contact pressure density in the neighborhood of the boundary $\Gamma_*$, $\mbox{\boldmath $\eta$}(s_*)$ is a point of the contour $\Gamma_*$ with the dimensionless natural coordinate $s_*$, while the kernel $M(t)$ is given by
$$
M(t)=\int\limits_0^{+\infty}\frac{\mathcal{L}(u)}{u}\cos ut\, du.
$$

\begin{figure}[h!]
%\vskip-1.0cm    
    \centering
    \includegraphics[scale=0.4]{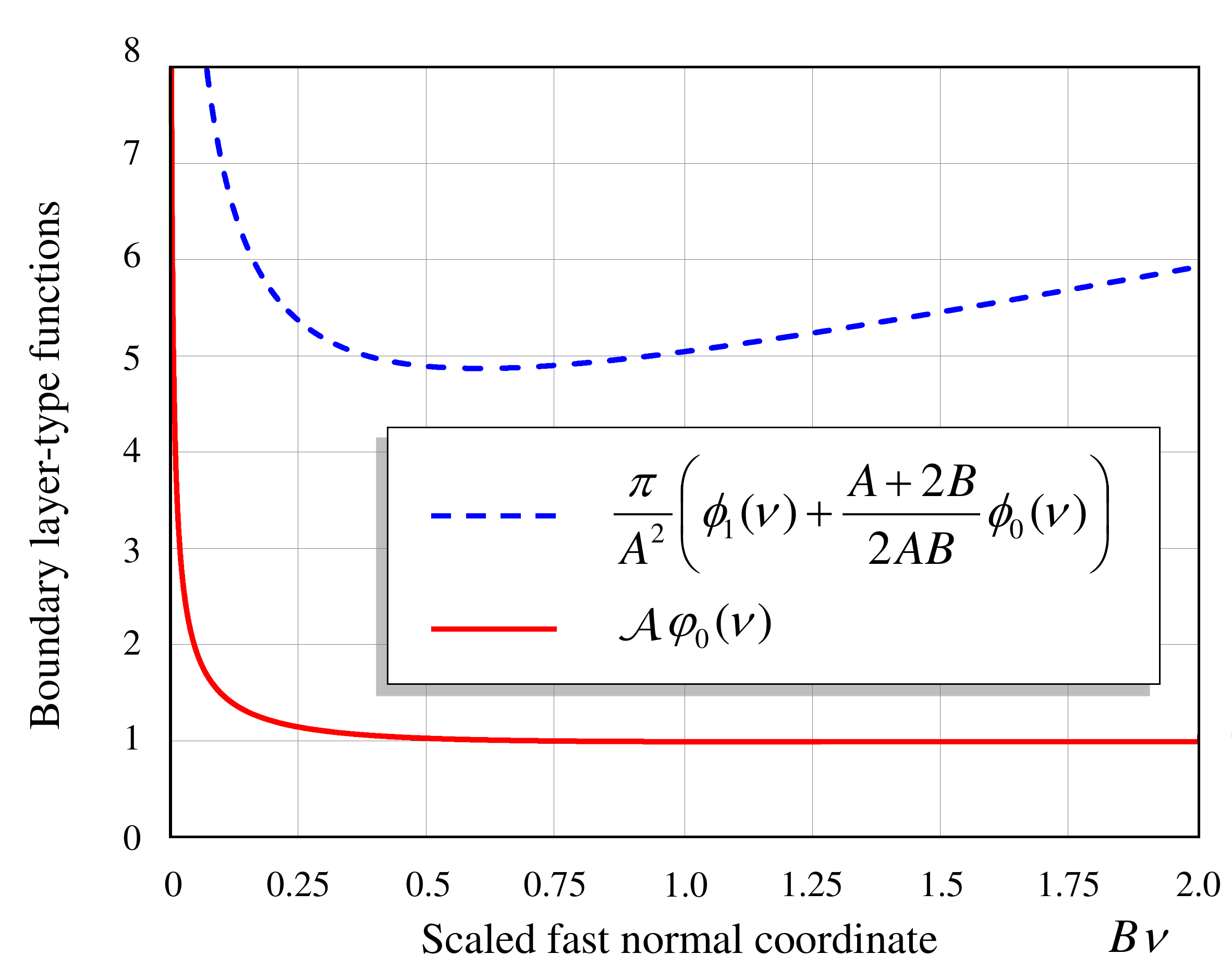}
%\vskip-3.5cm    
    \caption{Boundary-layer solutions normalized at the infinity. In the compressible case, the function (\ref{jkr1(3.8)}) is plotted for a typical value ${\mathcal{A}}B=0{.}5$. In the incompressible case, the normalized function (\ref{jkr1(4.10)}), with the regularity condition (\ref{jkr1(4.11)}) taken into account, for the isotropic material with $A=0{.}761310$ and $B=2{.}588024$ according to Aleksandrov's approximation~\cite{Alexandrov2003}.}
%\vskip-1.0cm        
    \label{phi01adhe}
\end{figure}

Making use of Aleksandrov's approximate factorization for the function $w^{-1}\mathcal{L}(w)$ suggested in \cite{Alexandrov1969-2} (see also \cite{ArgatovMishuris2015}, Section~2.3.3 ), we readily find
\begin{equation}
\tilde{q}^*_\varepsilon(s_*,\nu)=\frac{\pi\theta h_*}{h_*}
\bigl[\delta_0^*-\varphi^*\bigl(\mbox{\boldmath $\eta$}(s_*)\bigr)\bigr]
\varphi_0(\nu).
\label{jkr1(3.7)}
\end{equation}
Here, $\varphi_0(\nu)$ is an approximate solution to Eq.~(\ref{jkr1(3.6)}) with the right-hand side 1 such that 
\begin{equation}
\varphi_0(\nu) = \frac{1}{\mathcal{A}}{\,\rm erf}\sqrt{B\nu}+\frac{e^{-B\nu}}{\sqrt{\pi \mathcal{A}\nu}}
\label{jkr1(3.8)}
\end{equation}
with ${\rm erf}(x)$ being the error function. The behavior of the boundary-layer function $\varphi_0(\nu)$ is shown in Fig.~\ref{phi01adhe}.

From (\ref{jkr1(3.2)}), (\ref{jkr1(3.2b)}), (\ref{jkr1(3.7)}), and (\ref{jkr1(3.8)}), it follows that 
\begin{equation}
K_I(s)\simeq -\sqrt{\frac{2}{h\mathcal{A}}}\,\theta
\bigl[\delta_0-\varphi\bigl({\bf y}(s)\bigr)\bigr],
\label{jkr1(3.9)}
\end{equation}
where $\mathcal{A}$ is given by (\ref{jkr1(2.4A)}).

Now, in light of (\ref{jkr1(2.7)}) and (\ref{jkr1(3.1)}), Eq.~(\ref{jkr1(3.9)}) can be rewritten in the form of the boundary condition for Eq.~(\ref{jkr1(3.1)}) as
$$
p\bigl({\bf y}(s)\bigr)=-\sqrt{\frac{2\mathcal{A}\Delta\gamma}{\theta h}}\,A_{33},
$$
or, taking into account formula $(\ref{jkr1(2.4)})_1$, as follows:
\begin{equation}
p\bigl({\bf y}(s)\bigr)=-\sqrt{\frac{2A_{33}\Delta\gamma}{h}}.
\label{jkr1(3.10)}
\end{equation}

Thus, in the framework of the leading-order asymptotic model for a thin compressible elastic layer, the JKR-type adhesive contact theory imposes the requirement of a constant contact pressure round the periphery of the contact area. 

\subsection{Elliptic frictionless adhesive contact for a thin compressible layer}
\label{Sec3.3}

In this section, we consider the adhesive contact problem (\ref{jkr1(3.1)}), (\ref{jkr1(3.10)}) in the case (\ref{jkr1(2.1)}). By substituting (\ref{jkr1(2.1)}) into Eq.~(\ref{jkr1(3.1)}), we readily obtain 
\begin{equation}
p({\bf y})=\frac{A_{33}}{h}\biggl(\delta_0-\frac{y_1^2}{2R_1}-\frac{y_2^2}{2R_2}
\biggr),
\label{jkr1(3.11)}
\end{equation}
while Eq.~(\ref{jkr1(3.10)}) determines the contour $\Gamma$ by the equation
\begin{equation}
\delta_0-\frac{y_1^2}{2R_1}-\frac{y_2^2}{2R_2}=
-\sqrt{\frac{2h\Delta\gamma}{A_{33}}}.
\label{jkr1(3.12)}
\end{equation}

It is obvious from Eq.~(\ref{jkr1(3.12)}) that the contact area $\omega$ is elliptic. Let the major semi-axis and the eccentricity of the contour $\Gamma$ be denoted by $a$ and $e$, so that the minor semi-axis is given by $b=\sqrt{1-e^2}a$. By simple calculations (assuming that $R_1\geq R_2$) we find
$$
e=\sqrt{1-\frac{R_2}{R_1}},
$$
\begin{equation}
a^2=2R_1\biggl(\delta_0+\sqrt{\frac{2h\Delta\gamma}{A_{33}}}\biggr).
\label{jkr1(3.14)}
\end{equation}

Now, by substituting the expansion (\ref{jkr1(3.11)}) into the double integral on the left-hand side of Eq.~(\ref{jkr1(2.8)}) and integrating over the domain $\omega$ bounded by the ellipse (\ref{jkr1(3.12)}), we evaluate the contact force
\begin{equation}
F=\frac{\pi A_{33}}{2h}ab
\biggl(\delta_0+\sqrt{\frac{2h\Delta\gamma}{A_{33}}}\biggr).
\label{jkr1(3.15)}
\end{equation}

Further, Eq.~(\ref{jkr1(3.14)}) yields the punch's displacement
\begin{equation}
\delta_0=\frac{a^2}{2R_1}-\sqrt{\frac{2h\Delta\gamma}{A_{33}}},
\label{jkr1(3.16)}
\end{equation}
whereas Eq.~(\ref{jkr1(3.15)}), in light of (\ref{jkr1(3.16)}), can be represented in the form
\begin{equation}
F=\frac{\pi A_{33}}{2h}\sqrt{\frac{R_2}{R_1}}a^2
\biggl(\frac{a^2}{2R_1}-2\sqrt{\frac{2h\Delta\gamma}{A_{33}}}\biggr).
\label{jkr1(3.17)}
\end{equation}

In the axisymmetric case, when $R_2=R_1=R$ and $b=a$, Eqs.~(\ref{jkr1(3.16)}) and (\ref{jkr1(3.17)}) reduce to 
\begin{equation}
\delta_0=\frac{a^2}{2R}-\sqrt{\frac{2h\Delta\gamma}{A_{33}}},\quad
F=\frac{\pi A_{33}}{4Rh}a^4
-\sqrt{2}\pi a^2\sqrt{\frac{A_{33}\Delta\gamma}{h}}.
\label{jkr1(3.18)}
\end{equation}

Finally, in the isotropic case, we have
$$
A_{33}=\frac{E(1-\nu)}{(1+\nu)(1-2\nu)},
$$
and formulas (\ref{jkr1(3.18)}) coincide with the known solution obtained by Yang \cite{Yang2006}.

\section{Adhesive contact in the case of a thin incompressible layer}
\label{Sec4}

\subsection{Leading-order asymptotic model for a thin incompressible elastic layer bonded to a rigid base}
\label{Sec4.1}

First of all observe that the coefficient $\mathcal{M}_0$ in the expansion (\ref{jkr1(2.3)}) vanishes for an incompressible material, and in the main asymptotic term formula (\ref{jkr1(2.5)}) reduces to
\begin{equation}
\frac{1}{\pi^2 h\theta}\iint\limits_{\omega}
p({\bf y}^\prime)K({\bf y}-{\bf y}^\prime)\,d{\bf y}^\prime
\simeq 
-\mathcal{M}_1 h^3\Delta_y p({\bf y}),
\label{jkr1(4.1)}
\end{equation}
where the coefficient $\mathcal{M}_1$ is given by $(\ref{jkr1(2.4)})_2$.

Therefore, by replacing the left-hand side of Eq.~(\ref{jkr1(2.2)}) with the right-hand side of (\ref{jkr1(4.1)}), we arrive at the differential equation 
\begin{equation}
-\mathcal{M}_1 h^3\Delta_y p({\bf y})
=\delta_0-\varphi({\bf y}),\quad {\bf y}\in\omega,
\label{jkr1(4.2)}
\end{equation}
where $\omega$ is the unknown contact area. 

It is interesting to note that Eq.~(\ref{jkr1(4.2)}) requires two boundary conditions in order to determine the domain $\omega$. In the case of a fixed contact area (see \cite{ArgatovMishuris2015,Barber1990,Alexandrov2003}) even for a flat-ended punch (when the solution to the governing integral equation (\ref{jkr1(2.2)}) has a square-root singularity), the following boundary condition has been imposed: $p({\bf y})=0$ for ${\bf y}\in\Gamma$. In the case of non-adhesive unilateral contact (see \cite{ArgatovMishuris2015,Chadwick2002,AteshianLaiZhuMow1994-22}), the zero-pressure gradient boundary condition is additionally imposed. 

In the case of adhesive contact, the boundary conditions can be formulated based on the boundary-layer problem. 

\subsection{Leading-order boundary-layer solution in the incompressible case}
\label{Sec4.2}

Following \cite{ArgatovMishuris2015}, we rewrite Eq.~(\ref{jkr1(2.2)}) in the form 
\begin{equation}
\frac{h}{\pi^2 \theta}\iint\limits_{\omega}
p({\bf y}^\prime)K_*({\bf y}-{\bf y}^\prime)\,d{\bf y}^\prime
=\mathcal{W}_0({\bf y}),
\label{jkr1(4.4)}
\end{equation}
where the kernel is given by 
$$
K_*({\bf y})=\iint\limits_0^{{}\ \ \ +\infty}
\frac{\mathcal{L}(\lambda)}{\lambda^3}\cos\frac{\lambda_1 y_1}{h}\cos\frac{\lambda_2 y_2}{h}\,d\lambda_1 d\lambda_2.
$$

In the leading asymptotic order (see \cite{ArgatovMishuris2015}, Section~2.7.2), we have
\begin{equation}
\mathcal{W}_0({\bf y})\simeq 
h^3\mathcal{M}_1 p({\bf y}).
\label{jkr1(4.5)}
\end{equation}

Following the asymptotic procedure developed in \cite{ArgatovMishuris2015}, we derive from Eq.~(\ref{jkr1(4.4)}) the following boundary-layer equation:
\begin{equation}
\int\limits_0^{+\infty}
q^{**}(s_*,\nu^\prime)M_*(\nu^\prime-\nu)\,d\nu^\prime=\frac{\pi\theta}{h_*^3}
\bigl(C_0(s_*)+\varepsilon C_1^*(s_*)\bigr).
\label{jkr1(4.6)}
\end{equation}
Here, $q^{**}(s_*,\nu)=\varepsilon^2 q^{*}_\varepsilon(s_*,\nu)$ is the leading-order approximation for the contact pressure density in the neighborhood of the boundary $\Gamma_*$ of the contact area $\omega_*$ in the dimensionless coordinates (\ref{jkr1(3.3)}), $\nu$ is the fast normal coordinate (\ref{jkr1(3.5)}), and the kernel $M_*(t)$ is given by
$$
M_*(t)=\int\limits_0^{+\infty} \frac{\mathcal{L}(u)}{u^3}\cos ut\,du.
$$

The right-hand side of Eq.~(\ref{jkr1(4.6)}) has been obtained from the expansion
\begin{equation}
\mathcal{W}_0({\bf y})=\mathcal{W}_0\bigl({\bf y}(s)\bigr)
+\varepsilon h_*\nu\frac{\partial\mathcal{W}_0}{\partial n}
\bigl({\bf y}(s)\bigr)
+\frac{\varepsilon^2}{2} h_*^2\nu^2\frac{\partial^2\mathcal{W}_0}{\partial n^2}
\bigl({\bf y}(s)\bigr)+\ldots,
\label{jkr1(4.7)}
\end{equation}
so that
\begin{equation}
C_0(s_*)=\mathcal{W}_0\bigl({\bf y}(s)\bigr),\quad
C_1^*(s_*)=h_*\nu\frac{\partial\mathcal{W}_0}{\partial n}
\bigl({\bf y}(s)\bigr).
\label{jkr1(4.8)}
\end{equation}

Based on Aleksandrov's approximation \cite{Alexandrov2003} for the kernel function $\mathcal{L}(u)$ given by the formula 
$$
\tilde{\mathcal{L}}(u)=\frac{\pi}{(u^2+A^2)\sqrt{u^2+B^2}},
$$
where it is assumed that (see for details \cite{ArgatovMishuris2015}, Section~2.6.5)
\begin{equation}
\frac{1}{A^2 B}=\theta\mathcal{M}_1.
\label{jkr1(4.9)}
\end{equation}

Now, using the analytical solution obtained by Argatov and Mishuris \cite{ArgatovMishuris2015}, we readily get
\begin{equation}
q^{**}(s_*,\nu)=\frac{\pi\theta}{h_*^3}
\bigl[C_0(s_*)\phi_0(\nu)+\varepsilon C_1^*(s_*)\phi_1(\nu)\bigr],
\label{jkr1(4.10)}
\end{equation}
where
$$
\phi_0(t) = \frac{A^2B}{\pi}{\,\rm erf\,}(\sqrt{Bt})
+\frac{A\sqrt{B}e^{-Bt}}{2\pi^{3/2}t^{3/2}}(2At-1),
$$
$$
\phi_1(t) = \frac{A^2Bt}{\pi}{\,\rm erf\,}(\sqrt{Bt})
+\frac{e^{-Bt}}{4\pi^{3/2}\sqrt{B}t^{3/2}}
(4A^2 Bt^2-2A^2 t+A+2B).
$$

Observe that the functions $\phi_0(t)$ and $\phi_1(t)$ have singularities of the order $O(t^{-3/2})$ as $t\to 0^+$. That is why, the square root singularity of the approximation (\ref{jkr1(4.10)}) requires that
$$
-C_0(s_*)\frac{A\sqrt{B}}{2\pi^{3/2}}
+\varepsilon C_1^*(s_*)
\frac{(A+2B)}{4\pi^{3/2}\sqrt{B}}=0,
$$
from where it follows that 
\begin{equation}
C_0(s_*)=\varepsilon C_1^*(s_*)\frac{(A+2B)}{2AB}.
\label{jkr1(4.11)}
\end{equation}

Correspondingly, the coefficient at the asymptotic term of the order $O(t^{-1/2})$ as $t\to 0^+$, which according to (\ref{jkr1(2.6)}) is related to the SIF of the boundary-layer contact pressure density (\ref{jkr1(4.10)}), is given by 
\begin{equation}
\lim_{\nu\to 0^+}\nu^{1/2}q^{**}(s_*,\nu)
=\varepsilon C_1^*(s_*)
\frac{A\sqrt{B}}{\pi^{3/2}}\frac{\pi\theta}{h_*^3}.
\label{jkr1(4.12)}
\end{equation}

On the other hand, we have
\begin{equation}
q^*_\varepsilon(s_*,\nu)\sim
-\frac{K_I(s)}{\sqrt{2\pi\varepsilon h_*}\nu^{1/2}},\quad \nu\to 0^+.
\label{jkr1(4.13)}
\end{equation}

Therefore, from (\ref{jkr1(4.12)}) and (\ref{jkr1(4.13)}) it follows that 
\begin{equation}
\frac{K_I(s)}{\sqrt{2\pi\varepsilon h_*}}=
-\varepsilon^{-1} C_1^*(s_*)
\frac{A\sqrt{B}}{\pi^{1/2}}\frac{\theta}{h_*^3}.
\label{jkr1(4.14)}
\end{equation}

Now, taking into account the boundary condition (\ref{jkr1(2.7)}) and the relations (\ref{jkr1(4.9)}), we evaluate $C_1^*(s_*)$ from (\ref{jkr1(4.14)}) as follows:
\begin{equation}
C_1^*(s_*)=-\varepsilon^{1/2}h_*^{5/2}
\sqrt{2\mathcal{M}_1\Delta\gamma}.
\label{jkr1(4.15)}
\end{equation}

Hence, from (\ref{jkr1(4.8)}) and (\ref{jkr1(4.15)}), we obtain
\begin{equation}
\frac{\partial\mathcal{W}_0}{\partial n}
\bigl({\bf y}(s)\bigr)=-\varepsilon^{1/2}h_*^{3/2}
\sqrt{2\mathcal{M}_1\Delta\gamma}.
\label{jkr1(4.16)}
\end{equation}

Note also that for a transversely isotropic incompressible elastic layer, we have 
$$
\mathcal{M}_1=\frac{1}{3G^\prime},
$$
where $G^\prime$ is the out-of-plane shear elastic modulus. Fig.~\ref{phi01adhe} presents the behavior of the normalized boundary-layer solution (\ref{jkr1(4.10)}), with the condition (\ref{jkr1(4.11)}) taken into account, for the isotropic case according to Aleksandrov's approximation \cite{Alexandrov2003}.

Thus, in light of the fact that $C_0(s_*)$ is an order of magnitude smaller than $C_1^*(s_*)$ as $\varepsilon\to 0^+$ (see Eq.~(\ref{jkr1(4.11)})), based on Eqs.~(\ref{jkr1(4.2)}), (\ref{jkr1(4.5)}), (\ref{jkr1(4.16)}), and (\ref{jkr1(4.7)}), we derive the following leading-order asymptotic model for the adhesive frictionless unilateral contact for a thin incompressible elastic layer:
\begin{equation}
-\frac{h^3}{3G^\prime}\Delta_y p({\bf y})
=\delta_0-\varphi({\bf y}),\quad {\bf y}\in\omega,
\label{jkr1(4.18)}
\end{equation}
\begin{equation}
p({\bf y})=0,\quad 
\frac{\partial p}{\partial n}({\bf y})
=-\frac{1}{h^{3/2}}\sqrt{6G^\prime\Delta\gamma},\quad {\bf y}\in\Gamma,
\label{jkr1(4.19)}
\end{equation}
where $\partial/\partial n$ is the inward normal derivative.

Finally, we note that in the case of non-adhesive contact, when $\Delta\gamma=0$, the asymptotic model 
(\ref{jkr1(4.18)}), (\ref{jkr1(4.19)}) coincides with the leading-order asymptotic model of unilateral frictionless contact for a thin incompressible elastic layer developed in 
\cite{ArgatovMishuris2015,ArgatovMishuris2010e-6bip,ArgatovMishuris2011-22,Argatov2012multibody-22}.

\subsection{Axisymmetric frictionless adhesive contact in the incompressible case}
\label{Sec4.3}

In the axisymmetric case, the boundary-value problem (\ref{jkr1(4.18)}), (\ref{jkr1(4.19)}) reduces to
\begin{equation}
\frac{1}{r}\frac{d}{dr}\biggl(r\frac{dp(r)}{dr}\biggr)
=m(C r^2-\delta_0),\quad r<a,
\label{jkr1(4.20)}
\end{equation}
\begin{equation}
p(a)=0,\quad \frac{dp}{dr}(a)=\sqrt{2m\Delta\gamma},
\label{jkr1(4.21)}
\end{equation}
where we have introduced the auxiliary notation
\begin{equation}
m=\frac{3G^\prime}{h^3},\quad C=\frac{1}{2R}.
\label{jkr1(4.22)}
\end{equation}

Integrating Eq.~(\ref{jkr1(4.20)}) with the first boundary condition (\ref{jkr1(4.21)}) taken into account, we get
\begin{equation}
p(r)=\frac{m}{16}\bigl[C (r^2+a^2)-4\delta_0\bigr](r^2-a^2),
\label{jkr1(4.23)}
\end{equation}
while the second boundary condition (\ref{jkr1(4.21)}) implies that
$$
\frac{m}{4}(Ca^3-2\delta_0 a)=\sqrt{2m\Delta\gamma},
$$
from where it follows that 
\begin{equation}
\delta_0=\frac{Ca^2}{2}-\frac{2}{a}\sqrt{\frac{2\Delta\gamma}{m}}.
\label{jkr1(4.24)}
\end{equation}

Further, integrating the contact pressure density (\ref{jkr1(4.23)}), we readily obtain the contact force in the form 
$$
F=\frac{\pi m}{24}a^4(3\delta_0-Ca^2),
$$
or, in light of (\ref{jkr1(4.24)}), we finally get
\begin{equation}
F=\frac{\pi m}{48} a^4\biggl(Ca^2
-\frac{12}{a}\sqrt{\frac{2\Delta\gamma}{m}}\biggr).
\label{jkr1(4.25)}
\end{equation}

Apart from the notation (see Eqs.~(\ref{jkr1(4.22)})), formulas (\ref{jkr1(4.24)}) and (\ref{jkr1(4.25)}) coincide with the relations between the punch displacement $\delta_0$, the contact force $F$, and the contact radius $a$, previously derived by Yang \cite{Yang2002} in the isotropic case for a spherical punch.

\subsection{Adhesive contact with a slightly perturbed circular contact area}
\label{Sec4.4}

Let us assume that the punch shape function, $\varphi_\mu(r,\theta)$, is prescribed in polar coordinates $(r,\theta)$ as a sum
\begin{equation}
\varphi_\mu(r,\theta)=\varphi_0(r)+\mu\varphi_1(r,\theta),
\label{jkr1(5.1)}
\end{equation}
where $\varphi_0(r)$ is a monotonically increasing function of the radial coordinate $r$, such that $\varphi_0(0)=0$, the function $\mu\varphi_1(r,\theta)$ describes a small deviation of the punch surface from the axisymmetric shape, and $\mu$ is a small parameter.

In the polar coordinates, the asymptotic model (\ref{jkr1(4.18)}), (\ref{jkr1(4.19)}) takes the form 
\begin{equation}
\frac{\partial^2 p_\mu}{\partial r^2}+\frac{1}{r}\frac{\partial p_\mu}{\partial r}
+\frac{1}{r^2}\frac{\partial^2 p_\mu}{\partial \theta^2}
=m\bigl(\varphi_\mu(r,\theta)-\delta_0\bigr),\quad (r,\theta)\in\omega_\mu,
\label{jkr1(5.2)}
\end{equation}
\begin{equation}
p_\mu\bigr\vert_{\Gamma_\mu}=0,\quad
\frac{\partial p_\mu}{\partial n}\biggr\vert_{\Gamma_\mu}=
-\sqrt{2m\Delta\gamma},
\label{jkr1(5.3)}
\end{equation}
where the dimensional parameter $m$ is given by $(\ref{jkr1(4.22)})_1$.

In light of (\ref{jkr1(5.1)}), we assume that the boundary $\Gamma_\mu$ of the contact area $\omega_\mu$ can be described by the equation
$$
r=a+\mu h(\theta),\quad \theta\in[0,2\pi),
$$
where $\mu h(\theta)$ is a small variation of the circular contact area $\omega_0$ with the boundary $\Gamma_0$ described by the equation $r=a$, $\theta\in[0,2\pi)$.

At the same time, the domain $\omega_0$ is determined from the limit problem 
\begin{equation}
\frac{1}{r}\frac{d}{dr}\biggl(r\frac{dp_0(r)}{dr}\biggr)
=m\bigl(\varphi_0(r)-\delta_0\bigr),\quad r<a,
\label{jkr1(5.5)}
\end{equation}
\begin{equation}
p_0(a)=0,\quad \frac{dp_0}{dr}(a)=\sqrt{2m\Delta\gamma}.
\label{jkr1(5.6)}
\end{equation}

Following \cite{ArgatovMishuris2015} (see Section~8.3.2), the solution to the boundary-value problem (\ref{jkr1(5.5)}), (\ref{jkr1(5.6)}) can be obtained in the form 
$$
\frac{1}{m}p_0(r)=\frac{\delta_0}{4}(a^2-r^2)-\Theta_0(a,r),
$$
where we have introduced the notation 
$$
\Theta_0(a,r)=\int\limits_r^a\varphi_0(\rho)\rho\ln\frac{a}{\rho}\,d\rho
-\int\limits_0^r\varphi_0(\rho)\rho\ln\frac{r}{\rho}\,d\rho.
$$

Moreover, the punch displacement and the limit value ($\mu=0$) of the contact force, $F_0$, are given by
\begin{equation}
\delta_0=\frac{2}{a^2}\int\limits_0^r \varphi_0(\rho)\rho\,d\rho
-\frac{2}{a}\sqrt{\frac{2\Delta\gamma}{m}},
\label{jkr1(5.8)}
\end{equation}
\begin{equation}
\frac{1}{m}F_0=\frac{\pi}{4}\int\limits_0^a
\varphi_0(\rho)(2\rho^2-a^2)\rho\,d\rho
-\frac{\pi a^3}{4}\sqrt{\frac{2\Delta\gamma}{m}}.
\label{jkr1(5.9)}
\end{equation}

Following \cite{ArgatovMishuris2011p} (see also \cite{ArgatovMishuris2015}, Section~9.1.6), we express the solution to the perturbed adhesive contact problem (\ref{jkr1(5.2)}), (\ref{jkr1(5.3)}) as
$$
p_\mu(r,\theta)=p_0(r)+\mu p_1(r,\theta)+O(\mu^2),
$$
\begin{equation}
F_\mu=F_0+\mu F_1+O(\mu^2),
\label{jkr1(5.11)}
\end{equation}
assuming that the punch displacement $\delta_0$ is specified, while the contact force $F_\mu$ is unknown {\it a~priori\/}.

Applying the perturbation technique, we arrive at the following problem:
\begin{equation}
\frac{\partial^2 p_1}{\partial r^2}+\frac{1}{r}\frac{\partial p_1}{\partial r}
+\frac{1}{r^2}\frac{\partial^2 p_1}{\partial \theta^2}
=m\varphi_1(r,\theta),\quad r<a,
\label{jkr1(5.12)}
\end{equation}
\begin{equation}
p_1(a,\theta)=-h(\theta)\frac{dp_0}{dr}(a),
\label{jkr1(5.13)}
\end{equation}
\begin{equation}
\frac{\partial p_1}{\partial r}(a,\theta)=
-h(\theta)\frac{d^2 p_0}{dr^2}(a).
\label{jkr1(5.14)}
\end{equation}

Let us express the solution to the problem (\ref{jkr1(5.12)}), (\ref{jkr1(5.13)}) in the form 
\begin{equation}
p_1(r,\theta)=Y_1(r,\theta)+Y_0(r,\theta),
\label{jkr1(5.15)}
\end{equation}
where $Y_1(r,\theta)$ is the solution to Eq.~(\ref{jkr1(5.12)}) with the boundary condition $Y_1(a,\theta)=0$, while $Y_0(r,\theta)$ is the solution to the Dirichlet problem with the boundary condition (\ref{jkr1(5.13)}).

In view of the boundary condition $(\ref{jkr1(5.6)})_2$, the application of Poisson's integral yields 
\begin{equation}
Y_0(r,\theta)=-\frac{\sqrt{2m\Delta\gamma}}{2\pi}\int\limits_0^{2\pi}
\frac{(a^2-r^2)h(\theta^\prime)\,d\theta^\prime}{
a^2-2ra\cos(\theta-\theta^\prime)+r^2}.
\label{jkr1(5.16)}
\end{equation}

Therefore, the substitution of (\ref{jkr1(5.15)}) into Eq.~(\ref{jkr1(5.14)}) leads to the equation for determining the variation of the contact area
\begin{equation}
\frac{\partial Y_1}{\partial r}(a,\theta)
-\sqrt{2m\Delta\gamma}\bigl(\mathfrak{S}h\bigr)(\theta)
=-h(\theta)\frac{d^2 p_0}{dr^2}(a),
\label{jkr1(5.17)}
\end{equation}
where $\mathfrak{S}$ is the Steklov---Poincar\'e (Dirichlet-to-Neumann) operator.

Note also that, in light of Eqs.~(\ref{jkr1(5.5)}) and $(\ref{jkr1(5.6)})_2$, we have
\begin{equation}
\frac{d^2 p_0}{dr^2}(a)=m\bigl(\varphi_0(a)-\delta_0\bigr)
-\frac{\sqrt{2m\Delta\gamma}}{a}.
\label{jkr1(5.18)}
\end{equation}

Now, let the function $h(\theta)$ be represented by its Fourier series 
\begin{equation}
h(\theta)=\frac{a_0}{2}+\sum_{n=1}^\infty a_n\cos n\theta+b_n\sin n\theta,
\label{jkr1(5.19)}
\end{equation}
then according to (\ref{jkr1(5.16)}), we get
\begin{equation}
\bigl(\mathfrak{S}h\bigr)(\theta)=\frac{1}{a}
\sum_{n=1}^\infty n\bigl[ a_n\cos n\theta+b_n\sin n\theta\bigr].
\label{jkr1(5.20)}
\end{equation}

Thus, the substitution of (\ref{jkr1(5.18)})--(\ref{jkr1(5.20)}) into Eq.~(\ref{jkr1(5.17)}) allows to uniquely determine the Fourier coefficients $a_n$ and $b_n$ in terms of the contact radius $a$, provided Eq.~(\ref{jkr1(5.8)}) is taken into account. 

Finally, in view of the boundary condition $(\ref{jkr1(5.3)})_1$, it can be shown that the contact force (\ref{jkr1(5.11)}) is evaluated as
$$
F_\mu=F_0+\mu\int\limits_0^{2\pi}d\theta\int\limits_0^a
p_1(\rho,\theta)\rho\,d\rho+O(\mu^2),
$$
where $F_0$ and $p_1(r,\theta)$ are given by (\ref{jkr1(5.9)}) and (\ref{jkr1(5.15)}), respectively. 

\section{Discussion and conclusion}
\label{SecDC}

First of all, observe that the normalization of the boundary-layer contact density $q^{**}(s_*,\nu)$ with respect to the small parameter $\varepsilon$ depends on the normalization of the geometric parameters of the problem (see, e.g., formulas (\ref{jkr1(3.2b)})), which in turn govern the relative size of the contact area with respect to the elastic layer thickness. To simplify technical details of our analysis in the incompressible case, we did not consider this question in detail (we refer the reader to the book \cite{ArgatovMishuris2015}).

In the case of an incompressible elastic layer, the contact area beneath a punch shaped as an elliptic paraboloid is not elliptic, as it is the case for a compressible layer (see Section~\ref{Sec3.3}). That is why, solution of the adhesive contact problem (\ref{jkr1(4.18)}), (\ref{jkr1(4.19)}) in the case (\ref{jkr1(2.1)}) is of obvious special interest for further study.

Following \cite{ArgatovMishuris2015}, the constructed asymptotic models, which describe the adhesive contact between a thin elastic layer and a rigid punch, can be generalized to the case of contact between two thin elastic layers (compressible or incompressible). The case of contact between a compressible thin elastic layer with an incompressible one requires a special consideration. 

The leading-order asymptotic models of non-axisymmetric adhesive contact between a rigid punch and a thin elastic layer bonded to a rigid base constitute the main results of the present paper.

\bigskip

\section*{Acknowledgements}
The authors acknowledge support from the FP7 IRSES Marie Curie grant TAMER No~610547.
One of the authors (IIA) is also grateful to the DFG (German Science Foundation -- Deutsche Forschungsgemeinschaft) for financial support during his stay at the TU Berlin.

\end{document}